\documentclass{article}

\usepackage{amssymb,latexsym}

\usepackage[pdftex]{graphicx}

\usepackage{hyperref}

\begin{document}

\pagestyle{plain}

\newtheorem{theorem}{Theorem}[section]

\newtheorem{proposition}[theorem]{Proposition}

\newtheorem{lema}[theorem]{Lemma}

\newtheorem{corollary}[theorem]{Corollary}

\newtheorem{definition}[theorem]{Definition}

\newtheorem{remark}[theorem]{Remark}

\newtheorem{exempl}{Example}[section]

\newenvironment{example}{\begin{exempl}  \em}{\hfill $\square$

\end{exempl}}  \vspace{.5cm}

\renewcommand{\contentsname}{ }

\title{COLIN implies LIN for emergent algebras}

\author{Marius Buliga \\ \href{https://mbuliga.github.io}{mbuliga.github.io}}

\date{15.10.2021}

\maketitle

\section{Emergent algebras}

Emergent  algebras come from sub-riemannian geometry \cite{buligasub}, \cite{buligainf},  \cite{buligabraided}, where they have the name "dilation structure" or "dilatation structure", introduced in \cite{buligadil1}. 
At the end of this article there are references for emergent algebras and dilation structures. Mind that in all references I used consistently wrong the right quasigroup name, it should be left quasigroup everywhere. This will be fixed in the present exposition of the subject. 

Emergent algebras, first time introduced in \cite{buligairq},  are families of quasigroup operations indexed by a commutative group, which satisfy some algebraic relations and also topological (convergence and continuity) relations. Besides sub-riemannian geometry, they appear as a semantics of a family of graph-rewrite systems related to interaction combinators \cite{buligahis}, or lambda calculus \cite{buligaglc}. In \cite{buligaem} there is a lambda calculus version of emergent algebras. 

Emergent algebras need three ingredients: a commutative group $\Gamma$ with an invariant topological filter, 
an uniform topological space $X$ and a family of operations on $X$, indexed with the elements of $\Gamma$.  The operations satisfy algebraic axioms and also some continuity and convergence axioms which are the reason behind the denomination of these algebras as "emergent". 

Let's take the ingredients one by one.

\paragraph{The group of scalars.} This is a commutative group $\Gamma$ with an absolute, i.e. a topological filter on $\Gamma$ which is invariant with respect to the group multiplication. 

 The typical example is the group $\Gamma = (0,+\infty)$ with the  multiplication operation, where $0$, which is not an element of the group, but it is seen as the topological filter generated by the intervals $(0,a)$, with $a >0$. 
For any $x \in \Gamma$ and any interval $(0,a)$ the set of elemements from that interval multiplicated by $x$ (at left or at right, it does not matter because the group is commutative) is the interval $(0,ax)$, which implies that $0$, as a topological filter, is invariant with respect of the group operation. Such a topological filter is called an absolute. We shall use $0$ to write limits as $a \in \Gamma$ converges to $0$. 

Other examples of commutative groups with absolutes are $(\mathbb{Z},+)$ with $+\infty$ as absolute, or $\displaystyle \mathbb{C}^{*}$  with complex number multiplication and $0$ as absolute. 

\paragraph{The uniform space.} An uniform space $X$ is a nonempty set endowed with an uniformity. An uniformity is a topological filter $\mathcal{U}$ over $X \times X$ which satisfies the properties: 
\begin{enumerate}
  \item[(a)] the diagonal of $X$ $\displaystyle \left\{ (x,x) \mid x \in X \right\}$  is a subset of any $A \in \mathcal{U}$, 
  \item[(b)] for any $A \in \mathcal{U}$ we have $\displaystyle A^{-1} \in \mathcal{U}$ where $$\displaystyle A^{-1} = \left\{  (y,x) \mid (x,y) \in A \right\}$$
  \item[(c)] for any $A,B \in \mathcal{U}$ we have $AB \in \mathcal{U}$, where  $$\displaystyle AB = \left\{ (x,w) \mid \exists y \in X \mbox{ s.t. } (x,y) \in A \mbox{ and } (y,w) \in B \right\}$$
\end{enumerate}

The typical example of an uniform space is a metric space $(X,d)$, where the uniformity filter is generated by sets $A \subset X \times X$ defined for any $a > 0$ as
$$A = \left\{ (x,y) \mid d(x,y) < a \right\}$$
Then (a) is implied by $d(x,x) = 0$ for any $x \in X$, (b) is implied by $d(x,y) = d(y,x)$ for any $x,y \in X$ and (c) is implied by $d(x,z) \leq d(x,y) + d(y,z)$ for any $x,y,z \in X$. 

We shall use the absolute $0$ of the group of scalars $\Gamma$ with the uniformity on $X$ in statements which contain "$\lim{a \rightarrow 0} f(a,x) = g(x)$ uniformly with respect to x", for some functions $\displaystyle f: \Gamma \times X \rightarrow X$.

\paragraph{Definition 1.} 
An emergent algebra over a set $X$ is a family of idempotent left quasigroup operations over X, 
indexed by the group $\Gamma$, which satisfy the following algebraic and topological axioms (R1), (R2), (act) and (em). \\

For $\displaystyle a \in \Gamma$ let's denote the left quasigroup operations indexed by $a$ with $\displaystyle \circ_a$ and $\displaystyle \bullet_a$. Therefore 
$\displaystyle (X,\circ_a ,\bullet_a )$ is an idempotent left quasigroup: \\

\paragraph{(R1)} $\displaystyle x \circ_a x = x$ 
 
\paragraph{(R2)} $\displaystyle x \circ_a (x \bullet_a y) = x \bullet_a (x \circ_a y) = y$ \\

\paragraph{(act)} algebraic axioms which relate the operations: \\

      $\displaystyle x \circ_a ( x \circ_b y) = x \circ_{ab} y $
      
      $\displaystyle x \circ_1 y = y $
      
      $\displaystyle x \circ_{1/a} y = x \bullet_a y $\\

\noindent On $X$ we define the following operations which we will use later:\\ 

- the approximate difference \\  

$\displaystyle \Delta^{x}_{a} (y , z)  =  (x \circ_a y) \bullet_a (x \circ_a z) $\\ 

- the approximate sum \\  

$\displaystyle \Sigma^{x}_{a} (y, z) = x \bullet_{a} ( ( x \circ_{a} y) \circ_{a} z)$  \\ 

- the approximate inverse \\

$\displaystyle inv^{x}_{a} y = (x \circ_{a} y) \bullet_{a} x$ \\

\noindent Finally we have the topological axioms.

\paragraph{(em)} $X$ is an uniform space; as $a$ converges to $0$ \\

     $\displaystyle (x,y) \mapsto x \circ_a y$ converges uniformly to  $(x,y) \mapsto x$  

     $\displaystyle (x,y,z) \mapsto \Delta^{x}_{a} (y , z)$ converges uniformly to a function $(x,y,z) \mapsto \Delta^{x} (y , z)$

      $(\displaystyle x,y,z) \mapsto \Sigma^{x}_{a} (y , z)$ converges uniformly to a function $(x,y,z) \mapsto \Sigma^{x} (y , z)$ \\

\paragraph{Topological axioms, graph rewrites and "emergence".} The reason for topological axioms is to allow passage to the limit in any order. The axioms (R1), (R2), (act) can be reformulated into graph rewrites on trivalent graphs with typed nodes and numbered ports, which is the start of the chemlambda project, as described in \cite{buligahis}.  The axiom (em), translated into graph rewrites, says that some patterns (corresponding to the graphs of the approximate difference, sum and inverse) can be replaced with new nodes (for the limit operations) {\bf in any order}. From here "emerge" new algebraic identities, or equivalently new graph rewrites, involving the new operations. The same effect can be obtained, but without reference to uniform structures and limits, in the term rewriting formalism "em" \cite{buligaem}, which is an enhancement of lambda calculus. \\

\section{LIN, COLIN and SHUFFLE}

Here are some other relevant notations, along with the more usual names from the field of quasigroups.

\paragraph{Definition 2.} An emergent algebra is linear if it is {\bf left distributive}:\\ 

\paragraph{ (R3) or (LIN)} $\displaystyle x \circ_{a} (y \circ_{b} z) = (x \circ_{a} y) \circ_{b} (x \circ_{a} z)$ \\

\noindent An emergent algebra is co-linear if it is {\bf right distributive}:\\ 

\paragraph{(COLIN)} $\displaystyle (x \circ_{a} y) \circ_{b} z = (x \circ_{b} z) \circ_{a} (y \circ_{b} z)$ \\

\noindent An emergent algebra has the shuffle trick if it is {\bf medial}: \\

\paragraph{(SHUFFLE)} $\displaystyle (x \circ_{a} y) \circ_{b} (u \circ_{a} v) = (x \circ_{b} u) \circ_{a} (y \circ_{b} v)$ \\

An important example of a linear emergent algebra comes from conical groups. This is an important class of topological groups which can be seen as a generalizaton of vector spaces. 

\paragraph{Definition 3.} A conical group $\displaystyle (X,\cdot,e)$  is a topological group endowed with a multiplication by scalars operation $\displaystyle (a,x) \in \Gamma \times X \mapsto ax \in X$ such that the following algebraic relations are true: 

$\displaystyle  a (b x) = (ab) x$ \\

$\displaystyle  a (x \cdot y) = ax \cdot ay$ \\ 

$\displaystyle  a (x^{-1}) = ( ax)^{-1}$ \\

$ a e = e $\\
and there is a neighbourhood $U \subset X$ of the neutral element $e \in X$ and a neighbourhood $V \subset \Gamma$ of the filter $0$ such that 

$\displaystyle (a,x) \in V \times U \mapsto ax \in U$\\

$\displaystyle x  \in U \mapsto ax$ converges uniformly to $x \mapsto e$, as $a$ converges to $0$. \\

We have the following structure theorem, see theorem 10 \cite{buligadil1} or theorem 6.2 \cite{buligaem}. \\

\paragraph{Theorem 1.} For an emergent algebra $X$, fix an element $e \in X$ and define the operations: \\

(addition) $\displaystyle (x,y) \mapsto  x \cdot y = \Sigma^{e} (x , y)$  \\

(inverse) $\displaystyle x \mapsto x^{-1} = inv^{e} x = \Delta^{e} (x,e)$ \\

(scalar multiplication) $(a,x) \mapsto  a x = e \circ_{a} x$ \\

Then $\displaystyle (X,\cdot)$ is a conical group with neutral element $e$ which is conical. 

Conversely, for any conical group $X$ and for any $a \in \Gamma$ define: \\ 

$\displaystyle x \circ_{a} y =  x \cdot a ( x^{-1} \cdot y)$ \\ 

With this operation $X$ becomes a linear emergent algebra. \\ 

\paragraph{Curvature as deviation from linearity.} For an emergent algebra, the following term represents the deviation from linearity: \\

$\displaystyle LIN_{a,b}(x,y,z) =   y \bullet_{b} (x \bullet_{a} ((x \circ_{a} y) \circ_{b} (x \circ_{a} z)))$ \\

Indeed, (LIN) is equivalent with  $\displaystyle LIN_{a,b}(x,y,z) = z$. But in the realm of sub-riemanian geometry, more precisely for a dilation structure, this term is related to curvature, as explained in \cite{buligasub}, section 2.5 "Curvdimension and curvature". For an arbitrary element $x \in X$ and for $c \in \Gamma$ define \\

$\displaystyle R^{a}_{b,c}(x,u,v,w) = x \bullet_{a} LIN_{b,c} ((x \circ_{a} u, x \circ_{a} v, x \circ_{a} w) $\\

We know then, by Theorem 1, that $\displaystyle R^{a}_{b,c}(x,u,v,w)$ converges to $w$, as $a$ converges to $0$. In particular, for a riemannian manifold $X$, there is an associated emergent algebra given by \\

$\displaystyle x \circ_{a} exp_{x}(y) = exp_{x}(ay)$ \\

\noindent where $exp$ is the geodesic exponential. For this emergent algebra, we recognize the construction called Schild's ladder in the term  \\ 

$\displaystyle r^{a}_{x}(v,w) = log_{x} \left( R^{a}_{\frac{1}{2}, \frac{1}{2}}(x,x, exp_{x}(av), exp_{x}(aw)) \right)$ \\ 

The distance, measured in the tangent space at $x$, between $w$ and $\displaystyle r^{a}_{x}(v,w)$ is controlled by the absolute value of $\displaystyle a^{2} \langle R_{x}(v,w)v,w\rangle$, where $\displaystyle R_{x}$ is the Riemann curvature tensor at $x$.

\paragraph{Commutator as deviation from co-linearity.} We can measure the deviation from co-linearity with the term \\

$\displaystyle COLIN_{a,b}(x,y,z) =  (x \circ_{a} y) \bullet_{b} ((x \circ_{b} z) \circ_{a} (y \circ_{b} z))$ \\

For conical groups, the computation of this term shows that it is related with the commutator.  One can prove that for conical groups (COLIN) is equivalent with (SHUFFLE). 

The following theorem is a sort of Bruck-Murdoch-Toyoda theorem \cite{toyoda}, \cite{murdoch}, \cite{bruck}. 

\paragraph{Theorem 2.} Let $X$ be a conical group. The following are equivalent:
\begin{enumerate}
\item[-] $X$ is commutative

\item[-] The associated emergent algebra is medial, i.e. it satisfies (SHUFFLE) 

\item[-] The associated emergent alegbra is right distributive, i.e. satisfies (COLIN).   
\end{enumerate}

\paragraph{Proof.} In a conical group with neutral element $e$ we compute 

$$COLIN_{a,b^{-1}}(e,e \bullet_{a} y,z) = [y \circ_{b} e, z \circ_{a} e] z$$ 
where $\displaystyle [x,y] = x y x^{-1} y^{-1}$ is the commutator with respect to the group operation. But in conical groups the function $\displaystyle x \mapsto x \circ_{a} e$ is surjective for any $a \not = 1$. (COLIN) is therefore equivalent with commutativity of the group operation. We already know that (SHUFFLE) is equivalent with commutativity. \\

It is therefore natural to ask: 

\paragraph{Question.} Are there emergent algebras which satisfy (COLIN) but not (LIN)? \\

\section{An useful tool: geometric series}

The geometric series is

$$ \sum_{n=0}^{\infty} \varepsilon^{n}  = \frac{1}{1-\varepsilon}$$

for any $\varepsilon \in (0,1)$.

With the notation for dilations:

$$ \delta^{x}_{\varepsilon} y = x + \varepsilon (-x+y)$$

we can rephrase the geometric series as an existence result. Namely that the equation:

$$ \delta^{S}_{\varepsilon} 0 = x $$

has the solution

$$ S = \sum_{n=0}^{\infty} \left( \delta^{0}_{\varepsilon}\right)^{n} x$$

for any $\varepsilon \in (0,1)$.

The non-commutative version of this result is given in \cite{buligainf} proposition 8.4. 
In that article is proposed a non-commutative affine geometry, where usual affine spaces (over a vector space) are replaced with their non-commutative versions over a conical group.  The work uses dilation structures, which are metrical versions of emergent algebras. In the mentioned article  the proof uses the existence of a metric on the (non-commutative affine) space. We can do a simple proof without it.

Of course that the condition $\varepsilon \in (0,1)$ will be reformulated as $\varepsilon^{n} \rightarrow 0$ as $n \rightarrow \infty$.

Pick a base point $e$ , which will play the role of the $0$. Then we want to solve the equation in $S$

$$ \delta^{S}_{\varepsilon} e = x$$

for a given, arbitrary $x$ and for an $\varepsilon$ with the property that $\varepsilon^{n} \rightarrow 0$ as $n \rightarrow \infty$.

We want to prove that

(*) $$S = \sum_{n=0}^{\infty} \left( \delta^{e}_{\varepsilon}\right)^{n} x$$

where the sum is with respect to the non-commutative addition operation based at $e$. More precisely this operation is "emergent":

$$ \Sigma^{e} (v, w)  = \lim_{\varepsilon \rightarrow 0} \Sigma^{e}_{\varepsilon}(v,w)$$

where the approximate sum is

$$ \Sigma^{e}_{\varepsilon}(v,w) = \delta^{e}_{\varepsilon^{-1}} \delta^{\delta^{e}_{\varepsilon} v}_{\varepsilon} w$$

Recall that the dilations satisfy the (LIN) property:

(LIN) $$\delta^{e}_{\varepsilon} \delta^{x}_{\mu} y = \delta^{\delta^{e}_{\varepsilon} x}_{\mu} \delta^{e}_{\varepsilon} y$$

The conclusion (*) can be reformulated as: define $S_{0} = x$ and

$$ S_{n+1} = \Sigma^{e}(x, \delta^{e}_{\varepsilon} S_{n})$$

Then $S =\lim_{n \rightarrow \infty} S_{n}$.

But this is simple, due to the following identities coming from (LIN).

The first identity uses the fact that, once we defined the addition from dilations and a passage to the limit, then we can prove that dilations themselves express via addition. This gives the first identity:

$$ \Sigma^{e}_{\varepsilon}(v,w) = \Sigma^{e}(\delta^{v}_{\varepsilon} e, w)$$

The second identity is easier, just use (LIN), there is no passage to the limit involved.

$$ \Sigma^{e}_{\varepsilon}(v, \delta^{e}_{\varepsilon} w)  = \delta^{v}_{\varepsilon} w $$

With these identities, the recurrence relation of the non-commutative geometric series becomes:

$$ S_{n+1} = \Sigma^{e}(x, \delta^{e}_{\varepsilon} S_{n}) = \Sigma^{e}_{\varepsilon}(S, \delta^{e}_{\varepsilon} S_{n}) = \delta^{S}_{\varepsilon} S_{n}$$

therefore

$$ S_{n} = \left( \delta^{S}_{\varepsilon} \right)^{n} S_{0}$$

and the proof ends by recalling that $\varepsilon^{n} \rightarrow 0$ as $n \rightarrow \infty$, which implies

$$ \lim_{n \rightarrow \infty} S_{n} = \lim_{n \rightarrow \infty} \delta^{S}_{\varepsilon^{n}} S_{0} = \lim_{\varepsilon \rightarrow 0} \delta^{S}_{\varepsilon} S_{0} = S$$
The proof is new even in the commutative case.

Here is a non-commutative example. We are in the group $N$ of real $n \times n$ upper triangular matrices, with 1's on the diagonal. This is a subgroup of real linear group.

For any scalar $e \in (0,+\infty)$ we define a diagonal matrix

$$ \mathbf{e}_{ij} = e^{i} \delta_{ij}$$

(here $\delta_{ij} = 1$ if $i=j$, otherwise $\delta_{ij} = 0$; please don't make a confusion with the dilation, denoted also with the letter $\delta$, which is introduced in the following).

Conjugation with diagonal matrices is an automorphism of $N$, in particular for any $e \in (0,+\infty)$ and any $\mathbf{x} \in N$ we have

$$ \mathbf{e}^{-1} \mathbf{x} \mathbf{e} \in N$$

This allows us to define the dilation:

$$ \delta^{\mathbf{x}}_{e} \mathbf{y} = \mathbf{x} \mathbf{e}^{-1} \mathbf{x}^{-1} \mathbf{y} \mathbf{e}$$

We leave to the reader to prove the following:  
\begin{enumerate}
\item[(a)] this defines an emergent algebra over $N$

\item[(b)] that this is a linear emergent algebra.
\end{enumerate}

We can define now the geometric series with respect to the "addition" operation which is simply the matrix multiplication in the group $N$:  (here $I$ is the identity matrix, the neutral element of the group)

$$ \sum_{k=0}^{\infty} \delta_{e^{k}}^{I} \mathbf{x}$$

and by the previous post, it does converge if $e \in (0,1)$.

Let's compute the finite sums. First notice that

$$ \delta^{\mathbf{x}}_{e} \mathbf{y} = [\mathbf{x}, \mathbf{e}^{-1}]  \mathbf{e}^{-1} \mathbf{y} \mathbf{e}$$

where the square bracket denotes the commutator.

We then have:

$$ \sum_{k=0}^{m} \delta_{e^{k}}^{I} \mathbf{x} = \mathbf{x}\mathbf{e}^{-1} \mathbf{x} \mathbf{e} \mathbf{e}^{-2} \mathbf{x} \mathbf{e}^{2} ... \mathbf{e}^{-m} \mathbf{x} \mathbf{e}^{m}$$

or equivalently:

$$ \sum_{k=0}^{m} \delta_{e^{k}}^{I} \mathbf{x} = (\mathbf{x} \mathbf{e}^{-1})^{m} \mathbf{x} \mathbf{e}^{m}$$

All in all, our convergence of the geometric series result says that:

 if $e \in (0,1)$ and

$$ [\mathbf{y}, \mathbf{e}^{-1}] = x$$

then

$$ \mathbf{y} = \lim_{m \rightarrow \infty}  (\mathbf{x} \mathbf{e}^{-1})^{m} \mathbf{x} \mathbf{e}^{m}$$

\paragraph{Remark} that for $n=2$, i.e. for the case of $2 \times 2$ upper triangular matrices with 1's on the diagonal, the convergence result is the classical geometric series convergence. Compute by hand the partial sums in tis case to convince yourself that indeed the usual geometric (partial) sum appears in the $(1,2)$ position in the matrices.

\section{(COLIN) implies (LIN)}

In this section we give a negative answer to the question formulated in the first section. At first sight this is surprising, because (COLIN) and (LIN) are symmetric one to the other. Therefore it would be expected that, as we have examples of emergent algebras which satisfy (LIN) but not (COLIN), which are just the noncommutative conical groups, there should be examples of emergent algebras which satisfy (COLIN) but not (LIN). 

At closer examination we notice that the axiom (R2), satisfied in any emergent algebra, breaks the symmetry between (LIN) and (COLIN). However, in itself this remark is not sufficient to answer the question if there are emergent algebras which satisfy (COLIN) but not (LIN).

\paragraph{Theorem 3.} For an emergent algebra with the group $\Gamma = (0,\infty)$ the condition (COLIN) implies the condition (LIN)

\paragraph{Proof. }  \paragraph{Part 1.} Recall the (COLIN) condition: for any $a, b \in \Gamma$ and for any $x, y, z \in X$ we have

$$ (x \circ_{a} y) \circ_{b} z = (x \circ_{b} z) \circ_{a} (y \circ_{b} z)$$

Fix an element $e \in X$, otherwise arbitrary. (COLIN) is then equivalent with

$$ y \circ_{b} z = (e \circ_{b} z) \bullet_{a} ((e \circ_{b} y) \circ_{a} z)$$

If we replace $z$ with $e \circ_{a} z$ and then we use (R2) and some groupings of terms then we obtain the following relation equivalent with (COLIN):

$$ y \circ_{b} (e \circ_{a} z) = \Delta_{a}^{e}(e \circ_{b} z, E)$$

where $E$ is a relative operation, namely

$$ E = e \bullet_{a} ((e \circ_{a} y) \circ_{b} (e \circ_{a} z))$$

We pass now with $a$ to $0$ by using the topological axiom (em) and we obtain in the limit the following

$$ y \circ_{b} e = \Delta^{e}(e \circ_{b} z, E_{0})$$

where $E_{0}$ is the limit of $E$, therefore the infinitesimal dilation of coefficient $b$, based at $e$. We write it like this

$$ E_{0} = y \circ_{b}^{e} z$$

\paragraph{Part 2.} The relation (COLIN) passes to the infinitesimal level. Indeed, if we replace the operations with the infinitesimal operations (dilations) based at $e$ then (COLIN) remains true.  This is true because  for an arbitrary $c \in \Gamma$ we deduce from (COLIN) the relation

 $$ (x \circ_{a,c}^{e} y) \circ_{b,c}^{e} z = (x \circ_{b,c}^{e} z) \circ_{a,c}^{e} (y \circ_{b,c}^{e} z)$$

where

$$ x \circ_{a,c}^{e} y = e \bullet_{c} ((e \circ_{c} x) \circ_{a} (e \circ_{c} y))$$

We can then pass to the limit with $c$ to $0$ and we get the "infinitesimal" (COLIN) relation

$$ (x \circ_{a}^{e} y) \circ_{b} ^{e}z = (x \circ_{b}^{e}z) \circ_{a}^{e} (y \circ_{b}^{e} z)$$

But the infinitesimal emergent algebra based at $e$ satisfies (LIN) \cite{buligainf} proposition 7.11. Therefore now we know that it also satisfy (COLIN). As a consequence we get that it comes from a commutative conical group. We denote with a dot this conical group operation.

Because the group is conical and commutative it follows that

$$ (e \circ_{a} y) \cdot (e \circ_{b} y)  = e \circ_{a + b} y$$

via techniques explained in \cite{buligaem}. Therefore the conical group is a vector space, with group operation being vector addition and $e \circ_{a} x$ equal to the  scalar $a$ which multiplies $x$ (up to an arbitrary exponential).

\paragraph{Part 3. }The last two relations from Part 1 are then rewritten as

$$ y \circ_{b} e = (e \circ_{b} z^{-1})  \cdot y \cdot (e \circ_{b} (y^{-1}  \cdot z))$$

Commutativity of the group operation $\cdot$ and (LIN) for the infinitesimal level gives us the equivalent

$$ y \circ_{b} e = y \circ_{b}^{e} e = e \circ_{1-b} y$$

(where in the last equality we used Part 2 and we make a slight abuse of notation for the multiplication by the scalar $1-b$)

Now we come back to the initial (COLIN) and we remark that with the new knowledge we can rewrite it as

$$ e \circ_{1-a} (x \circ_{b} y) = (e \circ_{1-a} x ) \circ_{b} ( e \circ_{1-a} y)$$

which is equivalent with

$$ x \circ_{b} y = e \bullet_{1-a} ((e \circ_{1-a} x ) \circ_{b} ( e \circ_{1-a} y))$$

We pass to the limit with $a$ to $1$ this time and we obtain that

$$ x \circ_{b} y = x \circ_{b}^{e} y$$

therefore the emergent algebra is identical with the infinitesimal emergent algebra based at $e$. Therefore it satisfies (LIN). \\

We know more, actually, namely that (COLIN) is equivalent with (SHUFFLE). Indeed, we proved that (COLIN) implies that the emergent algebra is the one of a conical commutative group, or we already know that (SHUFFLE) is true if and only if we are in a conical commutative group.

\paragraph{Acknowledgements.} Thanks to D. Stanovsky, who points to his survey article \cite{stanovsky}, section 3, where some examples of non-medial but right-distributive quasigroups are given.

\end{document}